\documentclass[12pt,a4paper]{article}

\usepackage{latexsym}
\usepackage{amsfonts}
\usepackage{amssymb}

\renewcommand{\theequation}{\arabic{section}.\arabic{equation}}

\newcommand{\be}{\begin{equation}}
\newcommand{\ee}{\end{equation}}
\newcommand{\bea}{\begin{eqnarray}}
\newcommand{\eea}{\end{eqnarray}}

\def\G{{\cal G}}                                 %
\def\L{{\cal L}}                                 %
\def\R{{\cal R}}                                 %
\def\D{{\cal D}}                                 %
\def\H{{\cal H}}                                 %
\def\la{\langle}                                 %
\def\ra{\rangle}                                 %
\def\pa{\partial}                                %
\def\ad{\mathrm{ad}}                             %
\def\T{{\cal T}}                                 %
\def\Z{\mathbf{Z}}                               %
\def\C{\mathbf{C}}                               %
\def\U{{\cal U}}                                 %
\def\J{{\cal J}}                                 %
\def\W{{\cal W}}                                 %
\begin{document}

\vspace*{0.5cm}
\begin{center}
{\Large \bf Explicit description of twisted Wakimoto realizations \\
of  affine Lie algebras } 

\end{center}

\vspace{1.0cm}

\begin{center}

L. Feh\'er\footnote{Corresponding author, e-mail: lfeher@rmki.kfki.hu}
and B.G. Pusztai\\

\bigskip

{\em 
Department of Theoretical Physics,
 University of Szeged \\
Tisza Lajos krt 84-86, H-6720 Szeged, Hungary \\

}
\end{center}

\vspace{1.5cm}

\begin{abstract}

In a vertex algebraic framework, we present an explicit description 
of the twisted Wakimoto realizations of the affine Lie algebras 
in correspondence with an arbitrary finite order automorphism and a 
compatible integral gradation of a complex simple Lie algebra.
This yields generalized free field realizations of
the twisted and untwisted affine Lie algebras in any gradation. 
The free field form of the twisted Sugawara formula and examples 
are also exhibited.

\end{abstract}

\newpage

\section{Introduction}
\setcounter{equation}{0}

The Wakimoto type free field realizations \cite{Waki} of 
the affine Lie algebras \cite{Kac} proved very useful 
in many applications, for example in calculations 
of correlation functions in conformal field theories \cite{CFT}.
After Wakimoto's study of the $sl_2$ case,  
the construction was extended to all the untwisted affine Lie
algebras by Feigin and Frenkel \cite{FF}.
Among the large number of further investigations we here 
mention only  the papers \cite{Kur,ATY,FFR,KurTak} dealing 
with applications to WZNW correlation functions, 
and the article  \cite{dBF} 
that contains explicit formulas for the 
affine currents in terms of the free fields in the general case,  
which will be used in the present work.
Detailed reviews of the Wakimoto construction and its applications
can be found in \cite{BMP,Frenkel,CFT}.

The aim of the present work is to generalize the explicit formulas
of \cite{dBF} to the twisted affine Lie algebras.
The Wakimoto construction has been recently extended 
to the twisted case by  Szczesny \cite{Szczesny}.
Properly explained in the framework of vertex algebras, 
the basic idea is to induce the twist of the affine Lie 
algebra by twisting the 
free fields that enter the untwisted currents.
In this manner Wakimoto modules are constructed in \cite{Szczesny} 
for the twisted affine Lie algebras  using their
standard realization by means of a diagram automorphism of a simple 
Lie algebra.  
The description of the twisted currents in \cite{Szczesny} 
is not quite explicit, only 
the example $A_2^{(2)}$ is described  in a more or less explicit
manner. 
(The $sl_3$ case was considered in \cite{DGZ}, too.)

By combining the basic idea of \cite{Szczesny} with the explicit 
formulas of the untwisted Wakimoto currents in \cite{dBF}, we can  
easily derive explicit formulas 
that give the twisted affine currents
as composites of the twisted free fields 
in the general case.
In fact, we here perform the construction  
in correspondence with any finite order automorphism of
a simple Lie algebra, and present also the realization 
of the Sugawara formula in terms 
of the twisted free fields, which was  not considered in \cite{Szczesny}.

Let us recall that there is a natural correspondence between 
the $\Z$-gradations of the affine Lie algebras and
the finite order automorphisms of the simple Lie algebras.
Free field realizations compatible with a 
particular $\Z$-gradation of an affine Lie algebra should 
be useful, for example, to analyze WZNW orbifolds 
defined by the corresponding Lie algebra automorphism. 
They could also yield a convenient tool for the free field 
realization of the WZNW model under the boundary condition 
that the group-valued field belongs to a twisted loop group. 
The construction may be of interest in connection  with 
inner automorphisms of the simple Lie algebras, too.
This is illustrated for instance by the papers in \cite{Hara},
where a free field realization of $\hat{sl}_2$
compatible with the principal gradation is investigated.
See also \cite{KacTod,FerMan,BHO} for applications of 
affine Lie algebras twisted by inner automorphisms. 

The above remarks motivated us to describe  the twisted Wakimoto
construction in correspondence with arbitrary finite order
automorphisms of the simple Lie algebras. More precisely, 
in our `input data' the automorphism is complemented by
the choice of a compatible integral gradation of the simple
Lie algebra. In the untwisted case generalized Wakimoto
realizations associated with arbitrary 
integral gradations (or parabolic subalgebras) 
of the simple Lie algebras have been considered in \cite{FF,dBF}.
The most important case is that of the principal
gradation (Borel subalgebra) studied also in \cite{Szczesny}.

The paper is organized as follows.
In section 2 definitions and   
a useful lemma are presented.
Section 3 contains a recall of the explicit formulas of the
generalized `Wakimoto homomorphism' from \cite{dBF},
and a simple proof if its equivariance with respect to the actions of the
pertinent Lie algebra 
automorphism on the affine currents and on the free fields.
Section 4 is devoted to the twisted Wakimoto construction.
Our main result is the explicit formula of the 
`twisted Wakimoto homomorphism' given by Proposition 3.
We also discuss the realization of the Sugawara formula
in terms of the twisted free fields (see eq.~(\ref{4.30}) for the result)
and sketch a convenient way to describe the `input data'
that are needed in order to obtain examples.
Examples are presented in two out of the three appendices to the main text. 
Appendix B contains a description of the simplest $sl_2$ case,
which the reader may consult first to get a feeling 
about the twisted Wakimoto construction. 
In appendix C a generalized (non-principal) free field 
realization of $D_3^{(2)}$ 
is presented  as an illustration.

In addition to deriving new results, our aim also is to 
provide a self-contained description of the twisted Wakimoto 
construction, which may be useful for future studies of its applications.
Keeping this in mind,  we summarized 
relevant background information on (twisted modules of) vertex algebras 
briefly in subsection 4.1 as well as in a technical appendix.
The content of appendix A is well known to specialists in vertex algebras, 
but it is apparently less well known in the related physical literature dealing 
with twisted chiral algebras.

\section{Definitions and conventions}
\setcounter{equation}{0}

Let $\G$ be a complex simple Lie algebra with invariant 
scalar product $\la\ ,\ \ra$.
We shall use the generalized Wakimoto realizations \cite{FF,dBF}
of the untwisted affine Lie algebra based on $\G$ that can be associated  
with any $\Z$-gradation of $\G$.
Such a gradation is defined by the eigenvalues of $\ad_H$ for a diagonalizable
element $H\in\G$ according to  
\be
\G =\oplus_m\,  \G_m
\qquad
[\G_m, \G_n] \subset \G_{m+n}
\qquad\hbox{with}\quad
\G_m=\{\, X\in \G\,\vert\, [H, X] = mX\,\}.
\label{2.2}\ee 
Denoting  the subspaces of positive/negative grades  by $\G_\pm$,
this yields the decomposition 
\be 
\G=\G_- + \G_0 + \G_+.
\label{2.3}\ee 
In the general case ${\cal P}=(\G_0+\G_+)$
is a parabolic subalgebra and 
$\G_0$ is a reductive Lie algebra.
That is $\G_0$ decomposes into an Abelian 
factor, say  $\G_0^0$,
and simple factors, say $\G_0^i$ for $i>0$, 
\be 
{\cal G}_0 = \oplus_{i\geq 0}\, {\cal G}_0^i,
\label{2.5}\ee
where the factors are pairwise orthogonal with respect to $\la\ ,\ \ra$.
In the `principal case'  $\G_0$ is a Cartan subalgebra  and 
the parabolic subalgebra is a Borel subalgebra.

For the purposes of this paper we consider a pair $(\tau,H)$,
where $H$ is as above and $\tau$ is an automorphism of $\G$ that  
has finite order, denoted as $N$.
We assume that $\tau$ and the $\Z$-gradation are {\em compatible} 
in the sense that $\tau(H)=H$, which is clearly equivalent to 
\be
\tau(\G_m)= \G_m \qquad
\forall m\in {\bf Z}.
\label{2.7}\ee

\subsection{Convenient bases}

In our construction we shall use joint eigenbases of $(\tau, \ad_H)$.
We denote by $\{T_a\}$ and $\{T^a\}$ dual bases of $\G$,
$\la T_a, T^b\ra = \delta_a^b$,
such that the base elements are simultaneous 
eigenvectors of $\ad_H$ and $\tau$ with respective eigenvalues
designated as follows:
\be
[H, T_a] = h_a T_a,
\quad
[H,T^a]=h^a T^a,
\quad
\tau(T_a)= \omega_a T_a,
\quad
\tau(T^a)=\omega^a T^a.
\label{2.8}\ee
We have
\be
\omega^a = \exp(\frac{2\pi \mathrm{i}}{N} n^a),\quad
\omega_a = \exp(\frac{2\pi \mathrm{i}}{N} n_a),
\quad 
n^a, n_a \in \{0,1,\ldots, N-1\},
\label{2.9}\ee
\be 
h_a + h^a =0,
\qquad
\omega_a \omega^a =1
\qquad
\forall a=1,\ldots, \mathrm{dim}(\G).
\label{2.10}\ee

We shall also use dual bases 
$\{ L_\alpha\}\subset \G_-$ and $\{ U^\beta\}\subset \G_+$ 
that are assumed to satisfy the relations
$\la L_\alpha, U^\beta \ra =\delta_\alpha^\beta$ and
\be
[H, L_\alpha] = h_\alpha L_\alpha,
\quad
[H,U^\alpha]=h^\alpha U^\alpha,
\quad
\tau(L_\alpha)= \omega_\alpha L_\alpha,
\quad
\tau(U^\alpha)=\omega^\alpha U^\alpha.
\label{2.11}\ee
These notations are consistent by assuming that $T_a=L_\alpha$ for
$a=\alpha=1,\ldots, \mathrm{dim}(\G_-)$.

We need two kind of bases of $\G_0$. 
Firstly, $\{D^{i,\mu}\}$ and its dual $\{D_{\mu}^i\}$ stand 
for bases compatible with the decomposition $\G_0= \oplus_i \G_0^i$, that is,
$D^{i,\mu}$ and $D_{\mu}^i$ belong to $\G_0^i$. 
These base elements are not necessarily eigenvectors of $\tau$.
Secondly, $D^k$ and $D_{k}$ ($k=1,\ldots, \mathrm{dim}(\G_0)$)
stand for dual bases of $\G_0$ consisting of $\tau$-eigenvectors.
In this case the eigenvalues are denoted according to 
\be
\tau(D^k)= \omega_0^k D^k,
\quad
\tau(D_{k})=\omega_{0,k} D_{k},
\quad
\omega_0^k = \exp(\frac{2\pi \mathrm{i}}{N} n_0^k),\quad
\omega_{0,k} = \exp(\frac{2\pi \mathrm{i}}{N} n_{0,k})
\label{2.12}\ee
with $n_0^k, n_{0,k} \in \{0,1,\ldots, N-1\}$.
It is not difficult to see \cite{OV} that the $\tau$-eigenvectors in $\G_0$ 
can be chosen in such a way that each $D_k$ belongs either to
$\G_0^0$ or to a direct sum of a subset of the factors 
$\G_0^i$ (\ref{2.5}) formed  by pairwise isomorphic factors.

\subsection{Some polynomials}

Consider the simply connected Lie group $G_-$ whose Lie algebra is $\G_-$.
The general element $g_-\in G_-$ can be 
parametrized as\footnote{Summation over coinciding indices
in opposite position is understood throughout the paper.}
\be
g_- = e^q 
\quad\hbox{with}\quad 
q= q^\alpha L_\alpha \in \G_-,
\label{2.13}\ee
which provides a diffeomorphism between $G_-$ and $\G_-$.
The group $G_-$ is naturally represented on the Lie algebra $\G$,
induced by restricting the adjoint representation of a group $G$ 
associated with 
the Lie algebra $\G$.  
Using a symbolic notation that pretends that $G$ is a matrix group, 
this representation is given by the map 
$g_- \mapsto \R(g_-)\in \mathrm{End}(\G)$ for which 
\be
\R(g_-) X := g_- X g_-^{-1}
\qquad 
\forall X\in \G.
\label{2.14}\ee
We now define the polynomials $\T^a_{\ \,b}(q)$ in the 
variables $\{q^\alpha\}$ by 
\be
\R(e^{-q}) T^a = \T^a_{\,\ b}(q) T^b. 
\label{2.15}\ee   
Similarly, we set 
\be
\R(e^{-q}) U^\alpha = \U^\alpha_{\ \,b}(q) T^b,
\qquad
\R(e^{-q}) L_\alpha = \L_{\alpha,b}(q) T^b = \L_\alpha^{\ \,\beta}(q) L_\beta,
\label{2.16}\ee 
and 
\be
\R(e^{-q}) D^{i,\mu} = \D^{i,\mu}_{a}(q) T^a,
\qquad
\R(e^{-q}) D^k = \D^k_{\ \,a}(q) T^a.
\label{2.17}\ee
We also need the polynomials $\Phi_\alpha^\beta(q)$ 
defined by
\be
 \frac{\partial e^q}{\partial q^\alpha} e^{-q} 
= \Phi_\alpha^\beta(q) L_\beta ,
\label{2.18}\ee
and the polynomials  
$\Psi^{\alpha}_\beta(q)$
that form the inverse matrix, 
\be
\Phi_\alpha^\beta(q) \Psi_\beta^\gamma(q) = \delta_\alpha^\gamma.
\label{2.19}\ee 
Finally, we introduce the $\G$-valued polynomials $\Lambda_\alpha(q)$ by 
\be
\Lambda_\alpha(q)= 
\frac{\partial \Psi_\gamma^\lambda(q)}{\partial q^\alpha} 
\Phi_\lambda^\rho(q) \R(e^{-q}) [U^\gamma, L_\rho].
\label{2.20}\ee
It can be shown that $\Lambda_\alpha$ is actually $\G_-$-valued, and 
hence we can write 
\be
\Lambda_\alpha(q)  = 
\Lambda_{\alpha,b}(q) T^b =\Lambda_\alpha^{\ \beta}(q) L_\beta .
\label{2.21}\ee

We now establish certain homogeneity properties of the
above polynomials. For this let us observe that if 
$\sigma$ is any automorphism of $\G$ and $f(z)$ is any complex
power series, then $f(\ad_q)$ ($q\in \G_-$) yields a polynomial,
since $\ad_q$ is nilpotent, and this polynomial satisfies
\be
\sigma^{-1}\circ f(\ad_q)\circ \sigma 
= f(\ad_{\sigma^{-1}(q)}).
\label{2.22}\ee
Consider the transformation 
\be
q^\alpha \mapsto \omega^\alpha q^\alpha 
\quad\hbox{or equivalently}\quad
q=q^\alpha L_\alpha \mapsto \tau^{-1}(q).
\label{2.23}\ee

\medskip\noindent 
{\bf Lemma 1.} {\em The above-defined 
polynomials in the $\{q^\alpha\}$ obey the following relations:}
\begin{displaymath}
\U^\alpha_{\ \,b}(\tau^{-1}(q))=
\omega^\alpha \omega_b\U^\alpha_{\ \,b}(q),
\quad 
\L_{\alpha,b}(\tau^{-1}(q)) =
\omega_\alpha \omega_b \L_{\alpha,b}(q),
\quad 
\D^k_{\ \,a}(\tau^{-1}(q)) =
\omega_0^k \omega_a\D^k_{\ \,a}(q),
\end{displaymath}
\begin{displaymath}
\Phi^\alpha_\beta(\tau^{-1}(q)) = 
\omega^\alpha\omega_\beta \Phi^\alpha_\beta(q),\quad
\Psi^{\alpha}_\beta(\tau^{-1}(q)) =
\omega^\alpha \omega_\beta \Psi^{\alpha}_\beta(q),
\quad 
\Lambda_{\alpha,b}(\tau^{-1}(q))=
\omega_\alpha \omega_b \Lambda_{\alpha,b}(q). 
\end{displaymath}

\medskip\noindent 
{\em Proof.}  The relations in the first line are special cases of
\be
\T^a_{\,\ b}(\tau^{-1}(q)) =
\la \exp(-\ad_{\tau^{-1}(q)})(T^a), T_b\ra = 
\la \exp(-\ad_{q})(\tau (T^a)), \tau(T_b)\ra 
=\omega^a \omega_b \T^a_{\,\ b}(q),
\label{2.24}\ee
where we used (\ref{2.22}) and the $\tau$-invariance of $\la\ ,\ \ra$.
The other relations follow similarly since 
\be
\Phi^\alpha_\beta(q)=\la U^\alpha, f(\ad_q) L_\beta\ra,
\quad
\Psi_\beta^\alpha(q)= \la U^\alpha, f^{-1}(\ad_q) L_\beta \ra
\quad
\hbox{with}\quad
f(z)= \frac{e^z -1}{z},
\label{2.25}\ee 
where $f(z)$ and $f^{-1}(z)= \frac{1}{f(z)}$ 
are expanded into their Taylor series around $z=0$. {\em Q.E.D.}

We shall see that Lemma 1 implies an important equivariance property 
of the Wakimoto realizations under the assumption in (\ref{2.7}). 
For completeness, let us also present another homogeneity property, 
which is used in the construction of the Wakimoto realizations \cite{dBF}.  
Consider now the transformation
\be
q^\alpha \mapsto \lambda^{h^\alpha} q^\alpha 
\quad\hbox{or equivalently}\quad
q=q^\alpha L_\alpha \mapsto \sigma_\lambda^{-1}(q)
\quad\hbox{with}\quad 
\sigma_\lambda = e^{(\ln \lambda) \ad_H}
\label{2.26}\ee
for any $\lambda \in \C^*$. In the same way as Lemma 1, we obtain 

\medskip\noindent 
{\bf Lemma 2.} {\em The above polynomials, defined by using any basis of $\G$ consisting 
of eigenvectors of $\ad_H$, are homogeneous in the sense that}
\begin{displaymath}
\T^a_{\,\ b}(\sigma_\lambda^{-1}(q)) 
=\lambda^{h^a+h_b} \T^a_{\,\ b}(q), \qquad
\Phi^\alpha_\beta(\sigma_\lambda^{-1}(q)) = 
\lambda^{h^\alpha + h_\beta} \Phi^\alpha_\beta(q),
\end{displaymath}
\begin{displaymath}
\Psi^\alpha_\beta(\sigma_\lambda^{-1}(q)) = 
\lambda^{h^\alpha + h_\beta} \Psi^\alpha_\beta(q),\qquad
\Lambda_{\alpha,b}(\sigma_\lambda^{-1}(q))=
\lambda^{h_\alpha + h_b} \Lambda_{\alpha,b}(q).  
\end{displaymath}
\smallskip

Notice that 
$P(\{\lambda^{h^\alpha} q^\alpha\}) =
\lambda^n P(\{ q^\alpha\})$ with $n<0$ holds only for 
the identically zero polynomial  $P$, since $h^\alpha >0$ for any $\alpha$. 
This implies for example that $\Lambda_\alpha(q)$ varies 
in $\G_-$ as stated in (\ref{2.21}).

\section{Wakimoto homomorphism and its equivariance }
\setcounter{equation}{0}

In the first subsection we recall the explicit formula for the 
generalized Wakimoto realizations of the current algebras derived in 
\cite{dBF}. The construction relies only on the data $(\G,H)$,
the automorphism $\tau$ will be used in the second subsection.

\subsection{Recall of the generalized Wakimoto homomorphisms}

In correspondence with a basis $T_a$ of $\G$, 
let us consider the currents $J_a(z)$ that have the mode expansions 
\be
J_a(z) = \sum_{n\in \Z} J_a[n] z^{-n-1}
\label{3.1}\ee
and are subject to the commutation relations
\be
[J_a(z), J_b(w)] = \la [T_a,T_b], T^c\ra J_c(w) \delta(z,w) 
+ K \la T_a, T_b\ra \pa_w \delta(z,w).
\label{3.2}\ee
Here $K$ is a fixed (non-zero) complex number, $z$ and $w$ are 
formal variables and 
\be
\delta(z,w) = \sum_{n\in \Z} w^n z^{-n-1}.
\label{3.3}\ee
This formula encodes the affine Lie algebra spanned 
by the modes $J_a[n]$,
\be
[J_a[m], J_b[n]]= 
\la [T_a,T_b], T^c\ra J_c[m+n] +
 K \la T_a, T_b\ra m \delta_{m,-n}.
\label{3.4}\ee
By considering the vacuum Verma module of the affine Lie algebra 
in the usual manner, 
the fields $J_a(z)$ together with their derivatives and normal ordered 
products generate a vertex algebra (see e.g. \cite{Kacbeg,FB}), 
which we denote as $A(\G,K)$.

In association with the triangular decomposition in (\ref{2.3}),
the generalized Wakimoto realization of the current algebra relies on a
vertex algebra whose generating fields can
be labelled by an $H$-graded basis of $\G$ spanned by
\be
L_\alpha,\quad  U^\alpha,\quad  D^i_{\mu}  
\label{3.5}\ee
as introduced above. We denote the corresponding fields as 
\be
p_\alpha(z),\quad q^\alpha(z),\quad \jmath^i_{\mu}(z).
\label{3.6}\ee
They are postulated to have the mode expansions
\be
p_\alpha(z) = \sum_{n\in \Z} p_\alpha[n] z^{-n-1},
\quad
q^\alpha(z) = \sum_{n\in \Z} q^\alpha[n] z^{-n-1},
\quad
\jmath^i_{\mu}(z) = \sum_{n\in \Z} \jmath^i_{\mu}[n] z^{-n-1}
\label{3.7}\ee
together with the 
commutators\footnote{A conjugate pair $p_\alpha$, $q^\alpha$ subject to (\ref{3.8}) 
is often called a $\beta\gamma$-system in the literature.} 
\be
[ q^\alpha(z), p_\beta(w)]= \delta^\alpha_\beta \delta(z,w),
\label{3.8}\ee
and 
\be
[\jmath^i_{\mu}(z), \jmath_{\nu}^i(w)] = 
\la [D_{\mu}^i,D_{\nu}^i], D^{i,\theta} \ra \jmath^i_{\theta}(w) 
\delta(z,w) 
+ K_0^i \la D_{\mu}^i, D_{\nu}^i\ra \pa_w \delta(z,w).
\label{3.9}\ee
The numbers $K_0^i$ are determined in terms of $K$ by 
\be
2 K_0^0=
2K+{\vert \psi\vert^2} h^\vee= 2K_0^i+ {\vert \psi_i\vert^2} h_i^\vee 
\quad\hbox{for}\quad i>0.
\label{3.10}\ee
Here $h^\vee$ and $h_i^\vee$ are the dual Coxeter numbers of $\G$ and
the simple factors $\G_0^i$ in (\ref{2.5}), $\psi$ and $\psi_i$ stand for 
the respective  highest roots, whose length is defined by means of
$\langle\ ,\ \rangle$ and its restriction to $\G_0^i$.   
It is understood that all other commutators between the
basic fields vanish.
The construction of the vertex algebra 
generated by these fields, which we denote by $W(\G,K,H)$
($W$ for Wakimoto) 
is a well known matter \cite{FB}.

Let us recall that the product of two commuting fields is well-defined
in any vertex algebra, and in this case the product is associative.
Therefore we can uniquely associate a field $P(z)$ to any
polynomial $P(q)$ in the variables $q^\alpha$,
simply by replacing $q^\alpha$ by $q^\alpha(z)$ in the arguments of $P$
(this rule extends obviously to $\G$-valued polynomials as well).
If the fields $\phi(z)$ and $\psi(z)$ do not commute, then their normal
ordered product is the field $:\phi(z) \psi(z):$ given by
\be
:\phi(z) \psi(z): = \phi_+(z) \psi(z) + \psi(z) \phi_-(z),
\label{3.11}\ee
where 
\be
\phi(z)=\phi_-(z) + \phi_+(z),
\quad
\phi_-(z)=\sum_{n\in \Z_+} \phi[n] z^{-n-1} 
\quad\hbox{for}\quad 
\phi(z)=\sum_{n\in \Z} \phi[n] z^{-n-1}
\label{3.12}\ee
with $\Z_+=\{0,1,2,\ldots\}$.

After all this preparation, we can now define the 
`Wakimoto currents' $\J_a(z)$ to be the following fields 
in the vertex algebra $W(\G,K,H)$:
\be
\J_a(z) = - :p_\beta(z)
\left( \Psi_\alpha^\beta(z) \U^\alpha_{\ \,a}(z)\right):
+ \sum_i \jmath_{\mu}^i (z)\D^{i,\mu}_a(z) + 
\left(K \Phi_\alpha^\beta(z) \L_{\beta,a}(z) +
\Lambda_{\alpha,a}(z) \right) \pa_z q^\alpha(z).
\label{3.13}\ee 

\medskip\noindent
{\bf Proposition 1 \cite{FF,dBF}.} {\em The currents $\J_a(z)$ satisfy
\be
[\J_a(z), \J_b(w)] = \la [T_a,T_b], T^c\ra \J_c(w) \delta(z,w) 
+ K \la T_a, T_b\ra \pa_w \delta(z,w).
\label{3.14}\ee
Thus one obtains a vertex algebra homomorphism
$\W_H: A(\G,K) \rightarrow W(\G,K,H)$ by mapping the 
generating fields $J_a(z)$ according to}
\be
\W_H: J_a(z) \mapsto \J_a(z).
\label{3.15}\ee

It is clear that (\ref{3.15}) extends uniquely to a 
vertex algebra homomorphism 
since the currents $J_a(z)$ generate the 
vertex algebra $A(\G,K)$ in the sense
of the reconstruction theorem 
(Theorem 4.5 in \cite{Kacbeg} or Theorem 3.6.1 in \cite{FB}). 
We call $\W_H$ the `(generalized) Wakimoto homomorphism' 
associated with the parabolic 
structure  specified by the grading element $H$. 
This homomorphism was originally established 
by Feigin and Frenkel \cite{FF} using indirect arguments.
The explicit formula (\ref{3.13}) found in \cite{dBF} allows for a  
direct verification of (\ref{3.14}).

It is useful to collect the current components 
in the `$\G$-valued currents' $J(z)= J_a(z) T^a$,
$\J(z)=\J_a(z) T^a$ and to also introduce  the 
$\G_-$-valued, $\G_+$-valued and $\G_0$-valued 
fields $q(z)$, $p(z)$ and $\jmath(z)$ by   
\be
q(z)=q^\alpha(z) L_\alpha,
\quad
p(z)=p_\alpha(z) U^\alpha,
\quad
\jmath(z)= \sum_i \jmath_{\mu}^i(z) D^{i,\mu}.
\label{3.16}\ee
By definition, $J(z)$ is independent of the choice of the basis of $\G$, 
$q(z)$ is independent of the choice of the basis of $\G_-$ and so on. 
For example, if $D^k$ is an arbitrary basis of $\G_0$, then the equality 
\be
 \sum_i \jmath_{\mu}^i(z) D^{i,\mu}=\jmath(z)=\jmath_{k}(z) D^k 
\label{3.17}\ee
determines the fields $\jmath_{k}(z)$ as linear combinations of 
the $\jmath_{\mu}^i(z)$.
With this notation, the Wakimoto current $\J(z)$ takes the form   
\be
\J(z)= - :p_\beta(z)\left( \Psi_\alpha^\beta(z) 
(e^{-q} U^\alpha e^{q})(z) \right):
+ \jmath_k(z) (e^{-q} D^k e^{q})(z) + 
(K e^{-q} \frac{\pa e^{q}}{\pa q^\alpha}
 +\Lambda_{\alpha})(z) \pa_z q^\alpha(z),
\label{3.18}\ee 
where $e^{-q} \frac{\pa e^{q}}{\pa q^\alpha}$ is of course
a $\G_-$-valued polynomial.

\subsection{Equivariance of the Wakimoto homomorphism}

 From now on we assume that the grading element $H$ and the automorphism 
$\tau$ satisfy the compatibility condition (\ref{2.7}).
We can then choose the base elements in $\G$ to 
be simultaneously $\ad_H$ and  
$\tau$-eigenvectors with the eigenvalues denoted as in subsection 2.1. 
The map 
\be
\tau_A: J_a(z) \mapsto \omega_a J_a(z)
\label{3.19}\ee
respects the current algebra (\ref{3.2}) and thus it extends uniquely to 
an automorphism of the vertex algebra $A(\G,K)$.
This automorphism, also called $\tau_A$, has the same  order $N$ as $\tau$.
Similarly, we can define an automorphism $\tau_W$ of the vertex algebra
$W(\G,K,H)$ by extension of the map
\be
\tau_W: (p_\alpha(z), q^\alpha(z), \jmath_{k}(z)) \mapsto 
(\omega_\alpha p_\alpha(z), \omega^\alpha q^\alpha(z), 
\omega_{0,k} \jmath_{k}(z)).
\label{3.20}\ee 

\medskip\noindent 
{\bf Proposition 2.}
{\em If the automorphism $\tau$ and the gradation 
are compatible in the sense of (\ref{2.7}), then the Wakimoto homomorphism 
$\W_H: A(\G,K) \rightarrow W(\G,K,H)$ is $\tau$-equivariant:}
\be
\tau_W \circ \W_H= \W_H \circ \tau_A.
\label{3.21}\ee

\medskip\noindent 
{\em Proof.} 
Let us denote the `Wakimoto current' defined by the right hand side of 
(\ref{3.13}) as 
\be
\J_a(\{ p_\alpha(z)\}, \{ q^\alpha(z)\}, \{ \jmath_{k}(z)\})
\label{3.22}\ee
to emphasize that it is a composite expression formed out of the fields
listed as its arguments. 
Note that in terms of the homogeneous fields $\jmath_{k}(z)$ the second term 
in (\ref{3.13}) has the form 
\be
\jmath_{k} (z) \D^k_{\ \,a}(\{ q^\alpha(z)\}). 
\label{3.23}\ee
The statement of Proposition 2 is obviously equivalent to the relation  
\be
\J_a(\{ \omega_\alpha p_\alpha(z)\}, \{ \omega^\alpha q^\alpha(z)\}, 
\{ \omega_{0,k} \jmath_{k}(z)\}) =
\omega_a \J_a(\{ p_\alpha(z)\}, \{ q^\alpha(z)\}, \{ \jmath_{k}(z)\}).
\label{3.24}\ee
This holds as a result of the homogeneity 
properties of the polynomials that enter $\J_a(z)$ with respect 
to the transformation (\ref{2.23}), as described in Lemma 1.   
{\em Q.E.D.}
\medskip

The equivariance of the Wakimoto homomorphism has been proved in 
\cite{Szczesny}
 in the 
`principal special case' for which $\tau$ is a diagram automorphism and $H$ 
defines the principal grading of $\G$. 
By taking advantage of 
the explicit formula (\ref{3.13}) of \cite{dBF} we presented a simpler proof, 
which is valid in the general case for which $\tau(H)=H$.

\section{Wakimoto realizations of twisted current algebras}
\setcounter{equation}{0}

We below consider the situation for which Proposition 2 holds and
explain how the twisted modules of the vertex algebra 
$W(\G,K,H)$ then yield generalized Wakimoto realizations of the 
$\tau$-twisted current algebra.
The idea is the same that has been used by Szczesny in 
\cite{Szczesny}, where 
the special case for which $\tau$ is a diagram automorphism and 
$H$ defines the principal grading of $\G$ was considered. 
Our result is more general and our formulas are much more explicit.

\subsection{Useful identities for twisted modules of vertex algebras}

Let $N$ be a positive integer and $M$ a complex vector space.
An $N$-twisted field on $M$ is a formal series of the form 
\be
F(z) = \sum_{n\in \frac{1}{N} \Z} F[n] z^{-n-1},
\label{4.1}\ee
where $F[n] \in \mathrm{End}(M)$ and $F(z)v$ contains only finitely
many negative powers of the formal variable $z^{\frac{1}{N}}$ for any $v\in M$.

Let $\sigma$ be an automorphism of order $N$ of a vertex algebra $V$.
As a vector space, we identify $V$ with the corresponding space 
of fields (vertex operators) on $V$.
Let us recall (see \cite{FFR91,D,Li,BKT} and appendix A)
 that a $\sigma$-twisted module of $V$ is a linear map $R_\sigma$  
from the fields in $V$ to the $N$-twisted fields on some vector space $M$,
here denoted as
\be
R_\sigma: \phi(z) \mapsto \phi^\sigma(z),
\label{4.2}\ee 
for which   ${\mathrm{id}}_V$ maps to ${\mathrm{id}}_M$, 
\be
\sigma( \phi(z)) = e^{\frac{2\pi\mathrm{i}}{N} k} \phi(z)
\qquad (k\in \Z)
\label{4.3}\ee
implies that 
\be
\phi^\sigma(z) = \sum_{n\in \frac{k}{N}+ \Z} \phi^\sigma[n] z^{-n-1},
\label{4.4}\ee
and the so-called twisted Borcherds identity holds.

We shall use the following 
consequences of  the twisted Borcherds identity (see appendix A). 
Let $V^k$ stand for the homogeneous fields subject to  (\ref{4.3}),
$V^k=V^{k+N}$.
First, 
suppose that the fields $\phi(z)$ and $\psi(z)$ in $V$ have the commutator
\be
[\phi(z), \psi(w)]= \sum_{s=0}^{s_{\rm max}}  \chi_s(w) \pa_w^s \delta(z,w).
\label{4.5}\ee 
If $\phi\in V^k$, then this implies 
\be
[\phi^\sigma(z), \psi^\sigma(w)]= \sum_{s=0}^{s_{\rm max}}  
  \chi^\sigma_s(w) \pa_w^s \delta_k(z,w),
\label{4.6}\ee
where $\phi^\sigma$, $\psi^\sigma$, $\chi_s^\sigma$ are the twisted fields 
corresponding respectively to $\phi$, $\psi$, $\chi_s$ and 
\be
\delta_k(z,w)= w^{\frac{k}{N}} z^{-\frac{k}{N}}\delta(z,w),
\qquad
(\delta_k = \delta_{k+N}).
\label{4.7}\ee
Of course, if $\phi\in V^k$, $\psi\in V^l$, then $\chi_s\in V^{k+l}$ for any $s$.
Second, suppose again that $\phi\in V^k$, $\psi\in V$  and consider
the normal ordered product 
\be
\xi(z) = : \phi(z) \psi(z):  
\label{4.8}\ee
Then the twisted field corresponding to $\xi(z)$ is given by 
\be
\xi^\sigma(z) = :\phi^\sigma(z) \psi^\sigma(z): - 
\sum_{s=0}^{s_{\rm max}}   s! 
\left( \begin{array}{c} k/N \\ s+1 \end{array} \right) \chi_s^\sigma(z) z^{-s-1}
\label{4.9}\ee
where we use $\bigl( \begin{array}{c} \alpha \\ s+1 \end{array} \bigr)
=\frac{\alpha (\alpha -1) \cdots (\alpha -s)}{(s+1)!}$ 
for any $\alpha\in \C$, $s\in \Z_+$  and   
\be
:\phi^\sigma(z) \psi^\sigma(z): = \phi^\sigma_+(z) \psi^\sigma(z) + 
\psi^\sigma(z) \phi^\sigma_-(z),
\label{4.10}\ee
\be
\phi^\sigma_-(z)=\sum_{n\in \Z_+} \phi^\sigma[n+k/N] z^{-\frac{k}{N}-n-1}
\quad\hbox{with}\quad 
 k\in \{0,1,\ldots, N-1\}.
\label{4.11}\ee
In fact $:\phi^\sigma(z)\psi^\sigma(z):$  coincides with 
the `mode normal ordering' and $\xi^\sigma(z)$ with 
the `OPE normal ordering' 
of the twisted fields $\phi^\sigma(z)$, $\psi^\sigma(z)$ used 
in some references (see \cite{BHO,HW} and references therein).
If $[\phi(z), \psi(w)]=0$, then 
$\xi^\sigma(z)=\phi^\sigma(z)\psi^\sigma(z)$ as expected.
It can also be shown that for any field $\psi(z)$ in $V$ 
the twisted field corresponding to 
$\pa_z \psi(z)$ is given by $\pa_z \psi^\sigma(z)$.

\subsection{Application to twisted current algebras}

Let us suppose that we have a twisted module of the affine vertex algebra
$A(\G,K)$ corresponding to the automorphism $\tau_A$ (\ref{3.19}) and 
consider the twisted currents given by 
\be
R_{\tau_A}: J_a(z)\mapsto J_a^{\tau_A}(z) := J_a^\tau(z).
\label{4.12}\ee
These admit the mode expansions 
\be
J_a^\tau(z) = \sum_{n\in \frac{n_a}{N}+\Z} J_a^\tau[n] z^{-n-1}
\label{4.13}\ee
with the $n_a$ defined in (\ref{2.9}), and possess the commutators  
\be
[J_a^\tau(z), J_b^\tau(w)] = \la [T_a,T_b], T^c\ra J_c^\tau(w) \delta_{n_a}(z,w) 
+ K \la T_a, T_b\ra \pa_w \delta_{n_a}(z,w).
\label{4.14}\ee
according to  (\ref{3.2}), (\ref{4.6}). 
Equivalently, we obtain the `$\tau$-twisted current algebra' 
in the mode form 
\be
[J^\tau_a[m], J^\tau_b[n]]= 
\la [T_a,T_b], T^c\ra J^\tau_c[m+n] +
 K \la T_a, T_b\ra m \delta_{m,-n}.
\label{4.15}\ee
It is well known \cite{Kac} that this algebra is isomorphic to its untwisted 
analogue in (\ref{3.4}) if $\tau$ is an inner automorphism of $\G$,
and the isomorphism can be lifted to the corresponding vertex 
algebras \cite{Li}.   
Nevertheless, it is useful to have different realizations of the same algebra,
since they may have advantages from the viewpoint of various applications 
\cite{GO,Hara,KacTod,FerMan,BHO}.

Let us now consider a twisted module of the vertex algebra $W(\G,K,H)$ 
with respect to the automorphism $\tau_W$ (\ref{3.20}).
For the corresponding twisted fields now use the notation 
\be
R_{\tau_W}: (p_\alpha(z), q^\alpha(z), \jmath_{k}(z)) \mapsto 
(\tilde p_\alpha(z), \tilde q^\alpha(z), 
\tilde \jmath_{k}(z)).
\label{4.16}\ee  
The commutators of the basic twisted fields can be written as
\be
[ \tilde q^\alpha(z), \tilde p_\beta(w)]= \delta^\alpha_\beta \delta_{n^\alpha}(z,w),
\label{4.17}\ee
\be
[\tilde\jmath_{k}(z), \tilde\jmath_{l} (w)] = 
\la [D_k,D_l], D^m \ra \tilde\jmath_{m}(w) 
\delta_{n_{0,k}}(z,w) 
+ K_k \la D_k, D_l\ra \pa_w \delta_{n_{0,k}}(z,w),
\label{4.18}\ee
where the integers $n^\alpha$, $n_\beta$ and $n_{0,k}$ are defined in 
subsection 2.1.
We here use $\tau$-eigenvectors $D_k$ that belong either to
$\G_0^0$ or to a direct sum of a subset of the factors 
$\G_0^i$ (\ref{2.5}) formed 
by pairwise isomorphic factors. 
Correspondingly, $K_k$ 
is either $K_0^0$ or one of the $K_0^i$ in (\ref{3.10}).
The mode expansions of these twisted fields 
can be written as
\be
\tilde p_\alpha(z) = \sum_{n\in \frac{n_\alpha}{N}+\Z} \tilde p_\alpha[n] z^{-n-1},
\quad
\tilde q^\alpha(z) = \sum_{n\in \frac{n^\alpha}{N}+\Z} \tilde q^\alpha[n] z^{-n-1},
\quad
\tilde\jmath_{k}(z) = \sum_{n\in \frac{n_{0,k}}{N}+ \Z} \tilde\jmath_{k}[n] z^{-n-1}.
\label{4.19}\ee
One can of course present (\ref{4.17}), (\ref{4.18}) as an equivalent 
Lie algebra of the modes. 
For any polynomial $P$ in the complex variables $\{ q^\alpha\}$, note that 
\be
\tilde P(z) \equiv R_{\tau_W}(P(z)) 
\label{4.20}\ee
is obtained simply by replacing the $q^\alpha$ with the 
$\tilde q^\alpha(z)$ in the arguments of $P$. 
We need also the notation 
\be
\pa_\alpha P \equiv \frac{\pa P}{\pa q^\alpha}.
\label{4.21}\ee
Now we are ready to state the main result of this paper.

\medskip\noindent 
{\bf Proposition 3.}
{\em Suppose that the automorphism $\tau$ and the gradation of $\G$
are compatible (\ref{2.7}), and consider a twisted 
module $R_{\tau_W}$ of $W(\G,K,H)$ with respect to $\tau_W$ (\ref{3.20}).
Then one obtains a twisted module of $A(\G,K)$ with respect to $\tau_A$ 
(\ref{3.19}) by the composition
\be
R_{\tau_A} \equiv R_{\tau_W}\circ \W_H 
\label{4.22}\ee
with 
$\W_H: A(\G,K) \rightarrow W(\G,K,H)$ being the homomorphism described in Proposition 1.
The twisted currents $J_a^\tau(z) = R_{\tau_W}(\J_a(z))$ are found from (\ref{3.13})
explicitly  as 
\bea
&& J_a^\tau(z) = - :\tilde p_\beta(z)
\left( \tilde \Psi_\alpha^\beta(z) \tilde \U^\alpha_{\ \,a}(z)\right): 
+ \tilde \Theta_a(z) z^{-1} 
\nonumber\\
&& \qquad\qquad \qquad + \tilde\jmath_{k} (z)\tilde \D^k_{\ \,a}(z) + 
\left(K \tilde \Phi_\alpha^\beta(z) \tilde \L_{\beta,a}(z) +
\tilde \Lambda_{\alpha,a}(z) \right) \pa_z \tilde q^\alpha(z)
\label{4.23}\eea 
where the $\Theta_a$  are polynomials in the $\{ q^\alpha\}$ 
defined by  
$\Theta_a \equiv -\frac{n_\beta}{N} \pa_\beta (\Psi_\alpha^\beta \U^\alpha_{\ \,a})$.
}

\medskip\noindent 
{\em Proof.}
The first part of the statement relies on the $\tau$-equivariance of $\W_H$ (\ref{3.21}),
and is then obvious from the definition of twisted modules of vertex algebras
(see appendix A). 
Formula (\ref{4.23}) is derived by applying the identities 
collected in the preceding subsection, and using that for any polynomial $P$ 
in the variables $\{ q^\alpha\}$ one has
\be
[p_\beta(z), P(w)] = -(\pa_\beta P)(w) \delta(z,w)
\label{4.24}\ee
as a consequence of the Wick theorem. {\em Q.E.D.}

We now wish to describe     
the image of the `Sugawara field' 
\be
S(z)= \frac{1}{2y} : J_a(z) J^a(z) : =  \frac{1}{2y} \eta^{ab} : J_a(z) J_b(z) : 
\quad\hbox{with}\quad 
2y \equiv 2K+{\vert \psi\vert^2} h^\vee
\label{4.25}\ee
in a twisted Wakimoto module (\ref{4.22}) of $A(\G,K)$. 
To do this, we first quote from 
\cite{FF,dBF} that under the Wakimoto homomorphism 
${\cal S}(z)\equiv \W_H(S(z))$ has the form  
\be
{\cal S}(z)= \frac{1}{2y} : \J_a(z) \J^a(z) : = - : p_\alpha(z) \pa_z q^\alpha(z):  
+ \frac{1}{2y} \sum_{i\geq 0} : \jmath^i_{\mu}(z) \jmath^{i,\mu}(z): 
+ \frac{1}{y} \pa_z \jmath_{Q}(z)
\label{4.26}\ee  
with 
\be
\jmath_{Q}(z) \equiv \jmath^0_{\mu} (z) \langle D^{0,\mu}, Q\rangle 
\qquad\hbox{for}\qquad
Q= \frac{1}{2} [U^\alpha, L_\alpha].
\label{4.27}\ee
It is worth remarking that $Q=\rho_\G - \sum_{i>0} \rho_{\G_0^i}$,
where the Weyl vector $\rho_\G$ (respectively $\rho_{\G_0^i}$) 
corresponds to
half the sum of the positive roots of $\G$ (respectively $\G_0^i$) under 
the scalar product $\langle\ ,\ \rangle$, which implies that $Q\in \G_0^0$.
Denoting the inverse of $\eta_{0,kl}\equiv \la D_k, D_l\ra$ by
$\eta_0^{kl}$, we have
\be
\sum_{i\geq 0} : \jmath^i_{\mu}(z) \jmath^{i,\mu}(z): = \eta_0^{kl} 
: \jmath_{k}(z) \jmath_{l}(z): 
\label{4.28}\ee

As a direct consequence of (\ref{4.9}), 
for any twisted module of $A(\G,K)$ with respect to $\tau_A$, 
the image of $S(z)$ under 
$R_{\tau_A}$ (\ref{4.12}) is given\footnote{In an equivalent 
mode form, formula (\ref{4.29})
was first discovered in \cite{KP}, 
see also exercise 12.20 in \cite{Kac}. 
The local form (\ref{4.29}) can be found  
in \cite{BHO,HW} under the label 
`twisted affine-Sugawara construction'.}
as follows:
\be
S^\tau(z)= \frac{1}{2y} \eta^{ab} : J_a^\tau(z) J_b^\tau(z) :  
- \frac{z^{-1}}{2yN} \eta^{ab} n_a \la [T_a,T_b], T^c\ra J_c^\tau(z)  
+z^{-2} \frac{K}{4 y N^2} n_an^a  
\label{4.29}\ee
with the $n_a$ defined in (\ref{2.9}).
$S^\tau(z)$ obeys the same commutation relations as $S(z)$ since $\tau_A(S(z)) = S(z)$. 
In our case the twisted 
module of $A(\G,K)$ is given by (\ref{4.22}), and thus
we can determine $S^\tau(z)$ by applying $R_{\tau_W}$ to the 
right hand side of (4.26).  
By using (\ref{4.9}), we find that  
\bea
&& S^\tau(z)  = - : \tilde p_\alpha(z) \pa_z \tilde q^\alpha(z): + 
\frac{1}{2y}  :\tilde\jmath_{k}(z) \tilde \jmath^k(z): 
- \frac{z^{-1}}{2yN} \eta_0^{kl} n_{0,k} \la [D_k,D_l], D^m\ra \tilde\jmath_{m}(z) 
\nonumber \\
&& \qquad
+ \frac{1}{y} \pa_z \tilde \jmath_{Q}(z)
 +\frac{z^{-2}}{2N^2} \bigl( n_\alpha n^\alpha 
+ \frac{1}{2y} K_k n_{0,k} n^{0,k} \bigr).
\label{4.30}\eea
Recall that $K_k$ is defined below (\ref{4.18});
the $n_{0,k}$ label the eigenvalues of $\tau$ on $\G_0$ 
according to (\ref{2.12}); 
the $n_\alpha$ and the $n^\alpha$ similarly label the
$\tau$-eigenvalues on  $\G_-$ and on $\G_+$.   
Note also that $\tilde \jmath_{Q}(z)$
is expanded in integral powers of $z$ since  
$\tau_W (\jmath_{Q}(z)) = \jmath_Q(z)$
as a consequence of $\tau(Q)=Q$.

\subsection{How to construct examples?}

Since the automorphisms of finite order are fully under control 
\cite{Kac},
one can easily exhibit compatible pairs
$(\tau,H)$ together with their natural eigenbases.
Below we sketch the construction of these `input data'
of the twisted Wakimoto realizations  in rather general terms.

Let $\G$ be a simple Lie algebra with Cartan subalgebra $\H$, corresponding 
set of roots $\Delta$, and simple roots $\alpha_i$ for $i=1,\ldots, r$.
Any {\em inner} automorphism of order $N$ can be conjugated to have the
form  
\be
\tau= \exp( \frac{2\pi \mathrm{i}}{N} \ad_\theta) 
\ee
with some $\theta\in \H$ for which $\alpha_i(\theta)\in \Z$. Then a compatible 
integer gradation is obtained by taking $H$ to be any element from $\H$
for which $\alpha_i(H)\in \Z$. A joint eigenbasis for $(\tau,H)$ is 
provided by a Cartan-Weyl basis of $\G$, spanned by Cartan elements $H_{\alpha_i}$ 
and root vectors $E_\alpha$ for $\alpha\in \Delta$ normalized according to 
$\langle E_\alpha, E_{-\alpha}\rangle =1$.
Of course, the corresponding eigenvalues depend on the actual choice of $\tau$ and $H$.
A free field construction of the affine Lie algebra $\G^{(1)}$ 
in the  principal realization is obtained by 
taking $N$ to be the Coxeter number of $\G$ and setting both $\theta$ and $H$ equal to 
the element $I_0\in \H$ defined by 
\be
\alpha_i(I_0)= 1
\quad
\forall i=1,\ldots, r.
\label{I0}\ee
The simplest case of $\G=sl_2$ is 
described explicitly in appendix B.

Now recall \cite{Kac} that, up to conjugation, any {\em outer} 
automorphism of order $N$ can be written in the form
\be
\tau= \mu \circ \exp( \frac{2\pi \mathrm{i}}{N} \ad_\theta), 
\ee
where $\mu$ is induced by a non-trivial symmetry of the Dynkin diagram of $\G$,
$\theta$ belongs to the fixed point set of $\mu$ in $\H$ and 
 $\alpha_i(\theta)\in \Z$.
Then a compatible $\Z$-gradation of $\G$ can be specified by choosing any 
$H\in \H$ for which $\mu(H)=H$ and $\alpha_i(H)\in \Z$.
For simplicity of writing, let us assume that $\mu$ has order 2,
excluding only the cyclic diagram automorphism of $D_4$.
Denote by $\H^\pm$ and $\G^\pm$ the eigensubspaces of $\mu$ in $\H$ 
and in $\G$ in correspondence with the eigenvalues $\pm 1$. 
Recall that $\G^+$ is a simple Lie algebra and $\G^-$ is an irreducible
module of $\G^+$ in which the non-zero weights have multiplicity $1$.
Denote by $\Delta^+$ the roots of $(\H^+,\G^+)$ and by 
$\Delta^-$ the non-zero weights of $\H^+$ in $\G^-$. 
For $\lambda\in \Delta^\pm$ let $E^\pm_\lambda\in \G^\pm$ be an
eigenvector  of $\H^+$ normalized by $\langle E_\lambda^s, E_{-\lambda}^s\rangle =1$
for $s\in \{\pm\}$.
Choose also bases $\{ H_k^+\}_{k=1}^{r_+}$ of $\H^+$
and $\{ H_k^-\}_{k=1}^{r_-}$ of $\H^-$.
Now the elements
\be
E^+_\lambda\quad (\lambda\in \Delta^+), 
\quad
E^-_\lambda\quad (\lambda\in \Delta^-), 
\quad
H_k^+ \quad(1\leq k\leq r_+),
\quad
H_k^-\quad(1\leq k\leq  r_-)
\ee  
form a joint eigenbasis for any pair $(\tau,H)$ described above.
The subalgebra $\G_0$ equals $\H$ if we select $H:=I_0$ in (\ref{I0}).
A free field realization of the twisted affine Lie algebra $\G^{(2)}$ 
in the standard realization is obtained by taking $\tau:=\mu$.
This case has been studied in \cite{Szczesny}.
In our formalism it is equally easy to provide a free field construction 
of $\G^{(2)}$ in the principal realization.
It results from our general formulae by setting 
$H:=I_0$ together with $\theta:=I_0$ and $N$ twice the Coxeter
number of the algebra $\G^{(2)}$ (in the convention of \cite{Kac}).   

Our formalism contains also examples for which $\G_0$ is non-Abelian.
This is illustrated by a simple example displayed in appendix C. 

\section{Concluding remarks}

In this paper we described 
the generalized Wakimoto realizations of the twisted affine Lie algebras
in correspondence with arbitrary finite order automorphisms 
of the simple Lie algebras.
Our main result is the explicit formula given by Proposition 3, 
which relies on the equivariance of the formula of Proposition 1
valid in the untwisted case.
It is worth noting that these formulas admit classical limits
characterized by the absence of the terms 
represented by $\Lambda_\alpha(z)$ in (\ref{3.13}) and $\tilde \Theta_a(z)$ in (\ref{4.23}).
In fact, these terms are quantum corrections arising from  
the commutators of the basic free fields (the $\beta\gamma$-systems) that
contain the Planck constant if one uses suitable normalization.

In the classical WZNW model with twisted boundary condition, 
the classical twisted current algebra can be viewed 
as a consequence of a more fundamental quadratic `exchange algebra' of  
the group-valued twisted chiral WZNW field governed by a (non-unique) monodromy dependent 
dynamical $r$-matrix (see \cite{Sao} and references therein). 
The quantization of these quadratic Poisson algebras in terms of free 
fields is an interesting problem for the future.
It requires the free field realization of twisted chiral primary fields,
which should be also useful for further investigations of 
orbifolds and twisted versions of the WZNW model.

The method used in this paper is applicable in other cases as well,
for example to construct free field realizations of
twisted ${\cal W}$-algebras or twisted current algebras based on 
reductive Lie algebras that are a direct sum of identical simple factors \cite{BHS}.
Concerning the first case, recall \cite{dBT} that a ${\cal W}$-algebra can be associated  
with any $sl_2$-embedding into a simple Lie algebra $\G$.
If the $sl_2$ generators are invariant with respect to an automorphism $\tau$ of $\G$,
then $\tau$ lifts to an automorphism of the ${\cal W}$-algebra and 
the Miura maps underlying the free field realizations 
of the ${\cal W}$-algebra can be shown to be $\tau$-equivariant. 
In the second case the permutations of the identical factors give
rise to automorphisms of the current algebra, and its 
free field realization  obtained by applying the 
Wakimoto homomorphism (\ref{3.15}) factorwise is clearly equivariant with respect 
to the permutations.

\bigskip\bigskip
\noindent{\bf Acknowledgements.}
L.F. is indebted to J. de Boer for many useful discussions 
and for hospitality at ITFA during June 2002,
when a major part of this work was done.
He also wishes to thank  J. Balog, E. Ragoucy 
and J. Shiraishi for helpful remarks,
and B. Bakalov for sending him the manuscript \cite{BKT}.
This investigation was supported in part by the Hungarian 
Scientific Research Fund (OTKA) under T034170, T043159 and M036803.

\renewcommand{\theequation}{\arabic{section}.\arabic{equation}}
\renewcommand{\thesection}{\Alph{section}}
\setcounter{section}{0} 

\section{Vertex algebras and their twisted modules}
\setcounter{equation}{0}
\renewcommand{\theequation}{A.\arabic{equation}}

\def\Res{{\mathrm{Res}}}

For convenience, in this appendix we collect some background information on vertex algebras
and their twisted modules used in the main text.
For more details, see \cite{Kacbeg,FB,FFR91,D,Li,BKT} and references therein.

A vertex algebra is a quadruplet $(V,Y,v^0,T)$, where $V$ is a complex vector space
and $Y$ is a linear map from $V$ to the fields on $V$,
\be
Y: \phi \mapsto Y(\phi,z)= \sum_{n\in \Z} \phi_n z^{-n-1},
\quad \phi_n\in \mathrm{End}(V),
\quad \phi_n\psi =0\quad \forall  \phi, \psi\in V,\, n\gg 0,
\label{A.1}\ee
with $z$ being a formal variable.
By definition, 
the state-field correspondence $Y$, the vacuum vector $v^0\in V$, and the 
translation operator $T\in \mathrm{End}(V)$ obey the basic relations    
\be
Y(v^0,z)=\mathrm{id}_V, 
\quad
Tv^0 =0,\quad  [T, Y(\phi,z)]=\pa_z Y(\phi,z),
\label{A.2}\ee
\be
\phi_n v^0 = \delta_{n,-1} \phi
\quad
\hbox{for}\quad n\geq -1, \quad \forall \phi\in V, 
\label{A.3}\ee
as well as the Borcherds identity given by 
\bea
&&\Res_z \Bigl(  Y(\phi,z) Y(\psi,w) i_{z,w} F(z,w)
- Y(\psi,w) Y(\phi,z) i_{w,z} F(z,w)\Bigr) \nonumber\\
&& \qquad 
=\Res_{z-w} \Bigl(  Y(Y(\phi,z-w) \psi,w) i_{w,z-w} F(z,w)\Bigr)
\label{A.4}\eea 
for any $\phi,\psi\in V$ and $F(z,w)=  (z-w)^l z^m$ 
for any $l,m \in \Z$.
Here we have
\bea
&& i_{z,w} (z-w)^l=\sum_{n\in \Z_+} 
\left( \begin{array}{c} l \\ n \end{array} \right)
z^{l-n} (-w)^n,
\quad
i_{w,z} (z-w)^l=\sum_{n\in \Z_+} 
\left( \begin{array}{c} l \\ n \end{array} \right)
z^{n} (-w)^{l-n}, \nonumber\\
&& \qquad
i_{w,z-w} z^m =i_{w,z-w} ( w +(z-w))^m = \sum_{n\in \Z_+}
\left( \begin{array}{c} m \\ n \end{array} \right)
 w^{m-n} (z-w)^n.  
\label{A.5}\eea
Intuitively speaking, the Borcherds  identity (\ref{A.4}) 
(also called Cauchy-Jacobi identity) means that 
the `usual contour deformation' \cite{CFT} is applicable. 

Evaluating the Borcherds identity for $F(z,w)= (z-w)^l  z^m $ gives
\bea
&& \sum_{n \in \Z_+} (-1)^n 
\left( \begin{array}{c} l \\ n \end{array} \right)
\phi_{l+m -n} w^n  Y(\psi,w) 
+ Y(\psi,w) \sum_{n \in \Z_+} (-1)^{n+l +1} 
\left( \begin{array}{c} l \\ n \end{array} \right)
\phi_{m+n} w^{l-n} \nonumber\\
&& \qquad =\sum_{n\in \Z_+} \left( \begin{array}{c} m \\ n \end{array} \right)
w^{m-n} Y(\phi_{l+n} \psi,w).
\label{A.6}\eea
For $l=0$ (\ref{A.6}) simplifies to 
\be
\phi_m Y(\psi,w)- Y(\psi,w) \phi_m = \sum_{n\in \Z_+} 
\left( \begin{array}{c} m \\ n \end{array} \right)
w^{m-n} Y(\phi_n \psi,w),
\label{A.7}\ee
and the collection of these relations 
for all $m\in \Z$ is equivalent to the commutator formula
\be
[Y(\phi,z), Y(\psi,w)]= \sum_{n \in \Z_+} \frac{1}{n!} 
Y(\phi_n \psi, w) \pa_w^{n} \delta(z,w).
\label{A.8}\ee 
The sums on the right-hand-sides of the last two equations are actually finite;
$\delta(z,w)$ is defined in (\ref{3.3}).
In the special case $m=0$, $l=-1$ the Borcherds identity (\ref{A.6}) 
can be written as
\be
: Y(\phi,w) Y(\psi,w): = Y(\phi_{-1} \psi, w),
\label{A.9}\ee
and it can also be shown that $Y(T\phi,z) = \pa_z Y(\phi,z)$.

A vertex algebra homomorphism from $(V,Y,v^0,T)$ to 
$(\tilde V, \tilde Y, \tilde v^0, \tilde T)$ 
is defined to be a linear map $f: V \rightarrow \tilde V$ satisfying 
\be
f v^0 = \tilde v^0, 
\qquad
f \circ T = \tilde T \circ f,
\qquad
\tilde Y( f \phi, z)\circ  f  = f\circ Y(\phi, z) 
\quad \forall \phi\in V. 
\label{A.10}\ee
An isomorphism is an invertible homomorphism.  
An automorphism $\sigma$ of  $(V, Y, v^0, T)$ 
is thus defined by the properties 
\be
\sigma v^0 = v^0, 
\quad
[\sigma ,T] = 0,
\quad
Y( \sigma \phi, z) = \sigma \circ Y(\phi, z)\circ  \sigma^{-1}.  
\label{A.11}\ee
  
Let $\sigma$ be an automorphism of a vertex algebra 
$( V,  Y,  v^0,  T)$ that has finite order, say $N$.
Then 
\be
V=\oplus_{k=0}^{N-1} V^k
\quad\hbox{with}\quad 
V^k:= \{ \phi\in V\,\vert\, \sigma(\phi)= e^{\frac{2\pi \mathrm{i}k }{N}}\phi\}.
\label{A.12}\ee
A {\em twisted  $V$-module with respect to $\sigma$} 
is a vector space $M$ together with a linear map from $V$ to
the $N$-twisted fields on $M$,  
\be 
\phi \mapsto Y_M(\phi,z)=\sum_{n \in \frac{1}{N}\Z} \phi_n^M z^{-n-1},
\quad
\phi_n^M \in \mathrm{End}(M),
\quad \phi_n^M v =0\quad  \forall v\in M,\,\, n\gg 0,
\label{A.13}\ee
for which one has $Y_M(v^0,z) = \mathrm{id}_M$, 
\be
Y_M(\phi,z) = \sum_{n \in \frac{k}{N} + \Z}  \phi_n^M z^{-n-1}
\quad\hbox{if}\quad
\phi\in V^k,
\label{A.14}\ee
and the following twisted Borcherds identity:
\bea
&&\Res_z \Bigl( z^{\frac{k}{N}} Y_M(\phi,z) Y_M(\psi,w) i_{z,w}  F(z,w)
- Y_M(\psi,w) z^{\frac{k}{N}} Y_M(\phi,z) i_{w,z} F(z,w) \Bigr) \nonumber\\
&& \qquad 
=\Res_{z-w} \Bigl(  Y_M(Y(\phi,z-w) \psi,w) i_{w,z-w} z^{\frac{k}{N}} F(z,w)\Bigr)
\label{A.15}\eea 
$\forall \phi\in V^k$, $\psi\in V$ and $F(z,w)= (z-w)^l z^m $ 
$\forall l,m \in \Z$.
We here use  
\be
i_{w,z-w} (z^{\frac{k}{N}} z^m (z-w)^l)  =  
\sum_{n\in Z_+}  \left( \begin{array}{c} \frac{k}{N}+ m 
\\ n \end{array} \right) w^{\frac{k}{N}+m -n} (z-w)^{l+n}.
\label{A.16}\ee
The insertion of the factor $z^{\frac{k}{N}}$ ensures that the
`integrands' in the argument of $\Res$ contain only integral 
powers of the respective variables $z$ and $(z-w)$, 
and (\ref{A.15}) again can be thought of as the `contour deformation argument'. 
If $\sigma$ is the identity on $V$, then the 
twisted modules become just the (ordinary) modules of $V$. 

By evaluating the twisted Borcherds identity (\ref{A.15}) for $F(z,w)=(z-w)^l z^m$ 
one obtains 
\bea
&& \sum_{n \in \Z_+} (-1)^n 
\left( \begin{array}{c} l \\ n \end{array} \right)
\phi^M_{\frac{k}{N} + l+m -n} w^n  Y_M(\psi,w) 
+ Y_M(\psi,w) \sum_{n \in \Z_+} (-1)^{l +n+ 1} 
\left( \begin{array}{c} l \\ n \end{array} \right)
\phi^M_{\frac{k}{N}+m+n} w^{l-n} \nonumber\\
&& \qquad =\sum_{n\in Z_+} \left( \begin{array}{c} \frac{k}{N}+ m \\ n \end{array} \right)
w^{\frac{k}{N}+ m-n} Y_M(\phi_{l+n} \psi,w).
\label{A.17}\eea
For $l=0$ this simplifies to 
\be
\phi^M_{\frac{k}{N}+m} Y_M(\psi,w)- Y_M(\psi,w) \phi^M_{\frac{k}{N}+m}= \sum_{n\in \Z_+} 
\left( \begin{array}{c} \frac{k}{N}+m \\ n \end{array} \right)
w^{\frac{k}{N}+m-n} Y_M(\phi_n \psi,w),
\,\,\forall m\in \Z, 
\label{A.18}\ee
which is equivalent to the commutator formula
\be
[Y_M(\phi,z), Y_M(\psi,w)]= \sum_{n \in \Z_+} \frac{1}{n!} 
Y_M(\phi_n \psi, w) \pa_w^{n} \delta_k(z,w)
\quad 
\forall \phi\in V^k,\,\psi \in V.
\label{A.19}\ee
The sum is finite since $\phi_n \psi=0$ for $n\gg 0$; 
$\delta_k(z,w)$ appears in (\ref{4.7}).
The $m=0$, $l=-1$ special case of (\ref{A.17}) implies that 
\be
 Y_M(\phi_{-1} \psi, w)=: Y_M(\phi,w) Y_M(\psi,w): -
 \sum_{n \in \Z_+} \left( \begin{array}{c} \frac{k}{N} \\ n+1 \end{array} \right)
Y_M(\phi_n \psi,w) w^{-n-1}.
\label{A.20}\ee 
By definition, 
$:Y_M(\phi,w) Y_M(\psi,w):$ equals $w^{-\frac{k}{N}}$ times 
the expression 
on the left-hand-side of (\ref{A.17}) with $m=0$, $l=-1$.
It can also be shown that 
$Y_M(T\phi,z)= \pa_z Y_M(\phi,z)$ for any $\phi\in V$.

Upon changing the notation by setting 
\bea
&&Y(\phi,z)=\sum_{n\in \Z} \phi_n z^{-n-1} \equiv \phi(z) =
\sum_{n\in \Z} \phi[n] z^{-n-1} \nonumber\\
&&
Y_M(\phi, z)=\sum_{n\in \frac{1}{N} \Z} \phi^M_{n} z^{-n-1}
\equiv \phi^\sigma(z)=
\sum_{n\in \frac{1}{N}\Z} \phi^\sigma[n] z^{-n-1},
\qquad \forall \phi\in V, 
\label{A.21}\eea
the above-derived identities reproduce those mentioned in subsection 4.1.

\section{The simplest example: $\G=sl_2$}
\setcounter{equation}{0}
\renewcommand{\theequation}{B.\arabic{equation}}

The general features of the twisted Wakimoto construction
can be illustrated on the simplest $sl_2$ example as follows.
The untwisted current algebra based on $sl_2$ is generated by
\bea
&&[ J_{\sigma_3}(z), J_{\sigma_\pm}(w)]= \pm 2 J_{\sigma_\pm}(w)\delta(z,w),
\qquad
[J_{\sigma_3}(z), J_{\sigma_3}(w)]=2 K \pa_w \delta(z,w),\nonumber\\
&& [ J_{\sigma_+}(z), J_{\sigma_-}(w)] = J_{\sigma_3}(w) \delta(z,w) + 
K \pa_w \delta(z,w),
\label{E.1}\eea
where $J_{\sigma_a}(z)= \sum_{n\in \Z} J_{\sigma_a}[n] z^{-n-1}$
with the Pauli matrices
$\sigma_3= (E_{1,1}-E_{2,2})$, $\sigma_+=E_{1,2}$, $\sigma_-=E_{2,1}$.
Any automorphism of $sl_2$ of order $N$ is conjugate to 
\be
\tau= \exp( \frac{\pi \mathrm{i} }{N} \ad_{\sigma_3}),
\label{E.2}\ee
and  $\tau$ lifts to the automorphism $\tau_A$ of (\ref{E.1}) 
for which $\tau_A(J_{\sigma_a})= J_{\tau(\sigma_a)}$.
The corresponding twisted realization of the algebra (\ref{E.1}) is 
generated by the twisted currents $J^\tau_{\sigma_a}$ 
subject to the twisted commutation relations 
\bea
&&[ J^\tau_{\sigma_3}(z), J^\tau_{\sigma_\pm}(w)]= 
\pm 2 J^\tau_{\sigma_\pm}(w)\delta(z,w),
\qquad
[J^\tau_{\sigma_3}(z), J^\tau_{\sigma_3}(w)]=
2 K \pa_w \delta(z,w),\nonumber\\
&& [ J^\tau_{\sigma_+}(z), J^\tau_{\sigma_-}(w)] =
 J^\tau_{\sigma_3}(w) 
\bigl(w^\frac{1}{N} z^{-\frac{1}{N}} \delta(z,w)\bigr) + 
K \pa_w \bigl(w^\frac{1}{N} z^{-\frac{1}{N}} \delta(z,w)\bigr)
\label{E.3}\eea
and mode expansions 
$$
J^\tau_{\sigma_+}(z)= 
\sum_{n\in \frac{1}{N}+ \Z} J^\tau_{\sigma_+}[n] z^{-n-1},
\quad
J^\tau_{\sigma_-}(z)= 
\sum_{n\in \frac{N-1}{N}+ \Z} J^\tau_{\sigma_-}[n] z^{-n-1},
\quad
J^\tau_{\sigma_3}(z)= \sum_{n\in \Z} J^\tau_{\sigma_3}[n] z^{-n-1}.
$$

The free field  realization of (\ref{E.1}) due to \cite{Waki}
is given by the homomorphism $J_{\sigma_a}\mapsto \J_{\sigma_a}$ with 
the following `Wakimoto currents':
\bea
&& \J_{\sigma_-}(z)= -p(z),
\qquad 
\J_{\sigma_3}(z)=  \jmath(z) - 2 : p(z) q(z): \nonumber\\
&& \J_{\sigma_+}(z)=: p(z) (q^2)(z): - \jmath(z) q(z) + K \pa_z q(z),
\label{E.4}\eea
where $q$, $p$ and $\jmath$ satisfy the relations 
$[q(z), p(w)]=\delta(z,w)$ and 
$[\jmath(z), \jmath(w)]= 2(K+2)\pa_w \delta(z,w)$. 
On the free fields $\tau$ operates by the automorphism 
$\tau_W$ that maps $\jmath $ to itself and $(q,p)$ to 
$(\omega q, \omega^{-1} p)$ with $\omega=\exp(\frac{2\pi \mathrm{i}}{N})$.
The corresponding twisted free fields $\tilde q$ and $\tilde p$ satisfy 
\be
[\tilde q(z), \tilde p(w)] = 
w^\frac{1}{N} z^{-\frac{1}{N}} \delta(z,w),
\quad
\tilde q(z)= \sum_{n\in \frac{1}{N}+ \Z} \tilde q[n] z^{-n-1},
\quad
\tilde p(z)= \sum_{n\in \frac{N-1}{N}+ \Z} \tilde p[n] z^{-n-1}.
\label{E.5}\ee
It follows from Proposition 3 (eq.~(\ref{4.23})) that  
the twisted currents subject to (\ref{E.3})
can be realized in terms of the twisted free fields $\tilde q$,
$\tilde p$ and the (untwisted) $U(1)$ current $\jmath$ according to 
\bea
&& J^\tau_{\sigma_-}(z)= -\tilde p(z),
\qquad 
J^\tau_{\sigma_3}(z)=  \jmath(z) -
 2 : \tilde p(z)\tilde q(z): - 2z^{-1} \frac{N-1}{N} \nonumber\\
&& J^\tau_{\sigma_+}(z)=: \tilde p(z) (\tilde q^2)(z): - 
\jmath(z) \tilde q(z) + K \pa_z \tilde q(z) + 
2z^{-1} \frac{N-1}{N} \tilde q(z).
\label{E.6}\eea
The term $:\tilde p(z) (\tilde q^2)(z):+2z^{-1} \frac{N-1}{N} \tilde q(z)$
is the twisted analogue of $:p(z) (q^2)(z):$ as 
explained in the general case in Section 4.
See in particular eq.~(\ref{4.10}) for 
the definition of the normal ordering of the twisted 
free fields which is used here. 

The Wakimoto realizations of the twisted current algebra (\ref{E.3}) of $sl_2$
described in (\ref{E.6}) are different from the (similarly named) realization of 
$\hat{sl}_2$ in the principal gradation  presented in \cite{Hara}, 
which is based on different building blocks.
It would be interesting to find a connection between the two constructions.

\newpage
\section{An example with non-Abelian $\G_0$}
\setcounter{equation}{0}
\renewcommand{\theequation}{C.\arabic{equation}}

We below display a simple example for 
$\G:= sl_4$ such that $\G_0$ is non-Abelian. 
Since $A_3=D_3$, this yields a generalized free field realization 
of the affine Lie algebra $D_3^{(2)}$.

We consider the automorphism $\nu$ induced by transpose with 
respect to the `secondary diagonal':
\be
\nu: E_{i,j} \mapsto -E_{5-j, 5-i}
\label{B.1}\ee
for the usual elementary matrices $E_{i,j}\in gl_4$.
A compatible gradation is defined by 
\be
H:= \frac{1}{2}\left( E_{1,1} + E_{2,2} - E_{3,3} - E_{4,4}\right).
\label{B.2}\ee
Let $\H^+$ denote the diagonal matrices in $sl_4$ fixed by $\nu$,
and choose a basis of $sl_4$ consisting of weight vectors with 
respect to $\H^+$. An arbitrary element $D\in \H^+$ has the form
\be
D= d_1 E_{1,1} + d_2 E_{2,2} - d_2 E_{3,3} - d_1 E_{4,4}
\label{B.3}\ee
and we define the functional $e_i$ on $\H^+$ by $e_i(D)=d_i$.
A joint eigenbasis of $\nu$ and $\ad_H$ in $\G_-$ is furnished by 
the elements 
\be
E^-_{-2e_i} \equiv E_{5-i,i} \quad (i=1,2), \qquad
E^{\pm}_{-e_1 -e_2} \equiv E_{3,1} \mp E_{4,2}.
\label{B.4}\ee
Our notation indicates that $E_\lambda^\pm$ has weight $\lambda\neq 0$
with respect to $\H^+$ and $\nu$-eigenvalue $\pm 1$. 
With respect to $\langle X, Y\rangle \equiv {\mathrm{tr}}(XY)$ ($\forall X,Y\in \G$), 
the dual basis of $\G_+$ is spanned by 
\be
E^-_{2e_i} \equiv E_{i,5-i} \quad (i=1,2), \quad\hbox{and}\quad
\frac{1}{2}E^{\pm}_{e_1 +e_2}\quad \hbox{with}\quad 
 E^{\pm}_{e_1 +e_2}\equiv E_{1,3} \mp E_{2,4}.
\label{B.5}\ee
The base elements of $\G_0$ with non-zero weight are listed as   
\be
E^\pm_{e_1-e_2}\equiv E_{1,2}\mp E_{3,4},
\quad
E^\pm_{e_2-e_1}\equiv E_{2,1}\mp E_{4,3},
\label{B.6}\ee
and this is complemented by the basis of $\H\subset sl_4$ given by  
\be
H^\pm_{e_1-e_2}\equiv E_{1,1} - E_{2,2} \pm E_{3,3} \mp E_{4,4}
\label{B.7}\ee
together with $H$.
Note that $\G_0 = sl_2 \oplus sl_2 \oplus \mathrm{span}\{H\}$.

Let us use the expansions
\bea
&& q(z) =q^1(z) E^-_{-2e_1}+ q^2(z) E^-_{-2e_2} 
+ q^3(z) E^{+}_{-e_1 -e_2} +q^4(z) E^{-}_{-e_1 -e_2}, \nonumber\\
&& p(z) = p_1(z) E^-_{2e_1}+ p_2(z) E^-_{2e_2} 
+\frac{1}{2} p_3(z) E^{+}_{e_1 +e_2} +\frac{1}{2} p_4(z) E^{-}_{e_1+e_2},
\label{B.8}\eea
and the notation 
\be
\J_X(z)\equiv \J_a(z) \langle T^a, X\rangle
\quad \forall X\in \G,
\qquad
\jmath_{X}(z)\equiv \jmath_{k} \langle D^k, X\rangle
\quad \forall X\in \G_0.
\label{B.9}\ee
Dropping the variable $z$, the current 
components associated with $X\in \G_-$ are
\be
\J_{E^-_{-2e_i}}= - p_i\,\,\, (i=1,2),
\quad
\J_{E^{+}_{-e_1 -e_2}}= -p_3,
\quad
\J_{E^{-}_{-e_1 -e_2}}= -p_4.
\label{B.10}\ee 
For $X\in \G_0$ we find  
\bea
&&\J_{E^+_{e_1-e_2}} = \jmath_{E^+_{e_1-e_2}}
-2 :p_2 q^4: - :p_4 q^1:  \nonumber\\
&&\J_{E^-_{e_1-e_2}} = \jmath_{E^-_{e_1-e_2}}
-2 :p_2 q^3: + :p_3 q^1:  \nonumber\\
&&\J_{E^+_{e_2-e_1}} = \jmath_{E^+_{e_2-e_1}}
-2 :p_1 q^4: - :p_4 q^2:  \nonumber\\
&&\J_{E^-_{e_2-e_1}} = \jmath_{E^-_{e_2-e_1}}
+2 :p_1 q^3: - :p_3 q^2: \nonumber\\
&&\J_{H^+_{e_1-e_2}} = \jmath_{H^+_{e_1-e_2}}
-2 :p_1 q^1: +2 :p_2 q^2:  \nonumber\\
&&\J_{H^-_{e_1-e_2}} = \jmath_{H^-_{e_1-e_2}}
-2 :p_4 q^3: -2:p_3 q^4:  \nonumber\\
&&\J_{H} = \jmath_{H}
- :p_1 q^1: - :p_2 q^2: - :p_3 q^3: - :p_4 q^4:
\label{B.11}\eea 
For $X\in \G_+$ the explicit formulae are a bit longer. 
Straightforward calculation gives
\bea
&&\J_{E^-_{2e_1}} = -\frac{1}{2} q^1 \jmath_{ H^+_{e_1-e_2}} 
-q^1 \jmath_{ H} 
+ q^3 \jmath_{E^-_{e_1-e_2}}
- q^4 \jmath_{E^+_{e_1-e_2}}
 + : p_1 (q^1)^2 : 
\nonumber\\
&&\quad
+ : p_2 ((q^4)^2- (q^3)^2):
+ : p_3 (q^1 q^3) : 
+ : p_4 (q^1 q^4): +K \pa q^1,     
\label{B.12}\eea
\bea
&&\J_{E^-_{2e_2}} = \frac{1}{2} q^2 \jmath_{ H^+_{e_1-e_2}} -q^2 \jmath_{ H}
- q^3 \jmath_{E^-_{e_2-e_1}}
- q^4 \jmath_{E^+_{e_2-e_1}} + : p_2 (q^2)^2:
\nonumber\\
&&\quad
+ : p_1( (q^4)^2 - (q^3)^2): 
+ : p_3 (q^2 q^3) : 
+ : p_4 (q^2 q^4): +K \pa q^2,    
\label{B.13}\eea
\bea
&&\J_{E^+_{e_1+e_2}} = q^1 \jmath_{E^-_{e_2-e_1}} - q^2 \jmath_{E^-_{e_1-e_2}}
- 2 q^3 \jmath_{H} - q^4 \jmath_{H^-_{e_1-e_2}}
 + 2: p_1 (q^1 q^3): 
\nonumber\\
&&\quad
+ 2: p_2 (q^2 q^3):
+ : p_3 ( (q^3)^2 + (q^4)^2 - q^1 q^2): 
+ 2: p_4 (q^3 q^4): +2K \pa q^3,     
\label{B.14}\eea
\bea
&&\J_{E^-_{e_1+e_2}} =- q^1 \jmath_{E^+_{e_2-e_1}} - q^2 \jmath_{E^+_{e_1-e_2}}
- q^3 \jmath_{H^-_{e_1-e_2}} - 2 q^4 \jmath_{H}
 + 2: p_1 (q^1 q^4): 
\nonumber\\
&&\quad
+ 2: p_2 (q^2 q^4):
+ 2 : p_3 (q^3 q^4):
+ :p_4 ( (q^3)^2 + (q^4)^2 + q^1 q^2): 
 +2K \pa q^4.   
\label{B.15}\eea

The above formulae represent the Wakimoto current (\ref{3.13}) for 
$\G=sl_4$ with $H$ in (\ref{B.2}). 
The non-zero commutators  of the $\G_0$-valued current $\jmath(z)$
in the homogeneous basis of $\G_0$ are 
\bea
&&[ \jmath_{ H}(z), \jmath_{ H}(w)] =
(K+4) \pa_w \delta(z,w),\nonumber\\
&&[\jmath_{ H^s_{e_1-e_2}}(z), \jmath_{ H_{e_1-e_2}^s}(w)]=
4(K+2) \pa_w\delta(z,w), 
\nonumber\\
&& [ \jmath_{E^r_{e_1-e_2}}(z), \jmath_{E^s_{e_2-e_1}}(w) ]=
 \jmath_{ H_{e_1-e_2}^t}(w) \delta(z,w) + 2 (K+2) \delta_{r,s} 
\pa_w \delta(z,w), \nonumber\\
&& [ \jmath_{E^r_{\sigma (e_1-e_2)}}(z),
\jmath_{ H_{e_1-e_2}^s}(w)]= (-2 \sigma)
\jmath_{E^t_{\sigma (e_1-e_2)}}(w)\delta(z,w),
\quad \forall \sigma=\pm 1, 
\label{B.17}\eea
where $s,r \in \{\pm \}$, $t=+$ if $s=r$, and $t=-$ if $s\neq r$.  
Note that the action of $\nu$ on 
$\G_0=sl_2\oplus sl_2 \oplus \mathrm{span}\{H\}$ is given 
by the
transposition of the two $sl_2$ factors composed with an 
inner automorphism of $\G_0$.
On the vertex algebra generated by the components of $\jmath(z)$,
the corresponding automorphism is in fact 
a special case of the automorphisms  studied 
in relation to permutation orbifolds \cite{BHS,BDM}.
Thus one can construct the $\nu$-twisted modules of $W(\G,K,H)$ 
by applying the 
construction in \cite{BDM} to the current algebra
part and taking the usual twisted representations of the 
$\beta\gamma$-systems. 
In any such $\nu$-twisted module, the twisted Wakimoto 
currents $J^\nu_X$ ($X\in \G$) are obtained by substituting 
the twisted fields $\tilde \jmath_Y(z)$ ($Y\in \G_0$),
$\tilde q^\alpha$, $\tilde p_\alpha$ into 
(\ref{B.10}) -- (\ref{B.15}), and by also adding 
the correction terms defined  by the $\Theta_a$ in (\ref{4.23}). 
In our example, the non-zero correction terms associated 
with the base elements of $\G$ are
\be
\Theta_H=\frac{3}{2},
\quad
\Theta_{E^-_{2e_1}}=-\frac{3}{2} q^1,
\quad
\Theta_{E^-_{2e_2}}=-\frac{3}{2} q^2,
\quad
\Theta_{E^+_{e_1+e_2}}=-3 q^3,
\quad
\Theta_{E^-_{e_1+e_2}}=-3q^4,
\label{B.18}\ee
where we use $\Theta_X = \Theta_a \la X, T^a\ra$ for any $X\in \G$.

Finally, it is easy to describe 
the Wakimoto realizations of $D_3^{(2)}$
in correspondence with other automorphisms of $D_3=sl_4$, too.
For example, the diagram automorphism of 
$sl_4$ that leads to the standard realization (gradation) of $D_3^{(2)}$ 
can be written as 
\be
\mu = \nu \circ \exp( \pi \mathrm{i} \ad_{I_0}),
\qquad
I_0= \frac{3}{2} E_{1,1} + \frac{1}{2} E_{2,2} 
-\frac{1}{2}E_{3,3} -\frac{3}{2} E_{4,4},
\label{B.19}\ee
and the Coxeter automorphism of $sl_4$ associated with the
principal realization (gradation) of $D_3^{(2)}$ is given by
\be
\tau=\nu \circ \exp( \frac{4\pi \mathrm{i}}{3}  \ad_{I_0}).
\label{B.20}\ee
Since the above introduced basis of $sl_4$ is an eigenbasis 
of $\mu$ and $\tau$ as well, the twisted currents $J^\mu_a$
and $J^\tau_a$ are obtained by correspondingly twisting
the constituent fields in (\ref{B.10})--(\ref{B.15}) and adding
the appropriate correction terms $\Theta_a$ 
according to Proposition 3 (eq.~(\ref{4.23})).
These correction terms depend on the automorphism through
the integers $n_\alpha$ that label the 
eigenvalues of the automorphism on $\G_-$.

\end{document}